%%%%%%%This is a AMSTex file%%%%%%%%%%%%%%%%%

\input amstex
\documentstyle{amsppt}\magnification=1200
%\magnification=\magstep1
\NoBlackBoxes
\def\id{\operatorname{id}}
\def\Id{\operatorname{Id}}

\topmatter
\title 
Geodesic Conjugacy in two-step Nilmanifolds 
\endtitle
\rightheadtext{Geodesic Conjugacy in two-step nilmanifolds}
\author Carolyn Gordon and
Yiping Mao\endauthor
\address Department of Mathematics, Dartmouth College, 
         Hanover, NH 03755 \endaddress
\email  carolyn.s.gordon\@dartmouth.edu \endemail
\address Department of Mathematics, Texas Tech University,
         Lubbock, TX 79409\endaddress
\email ymao\@ttmath.ttu.edu \endemail
 \endkeywords
\thanks
The first author was supported in part by a grant from the National Science
Foundation.  Research at MSRI is supported
by NSF grant \#DMS 9022140.
\endthanks
\abstract
Two Riemannian manifolds are said to have $C^k$-conjugate geodesic flows if there exist an $C^k$ diffeomorphism between their unit tangent bundles which
intertwines the geodesic flows. We obtain a number of rigidity results for the geodesic flows on compact 2-step Riemannian nilmanifolds: For generic 2-step
nilmanifolds the geodesic flow is $C^2$ rigid. For special classes of 2-step nilmanifolds, we show that the geodesic flow is $C^0$ or $C^2$ rigid. In particular, there exist continuous families of 2-step nilmanifolds whose Laplacians are isospectral but whose geodesic flows are not $C^0$ conjugate.
\endabstract
\endtopmatter

\document

\heading  Introduction \endheading

Two Riemannian manifolds $(M, g)$ and $(N,h)$ are said to have 
$C^{k}$-{\it conjugate geodesic flows} if 
there is a $C^{k}$ diffeomorphism $F: S(M, g) \rightarrow S(N, h)$ which intertwines
the geodesic flows on $S(M, g)$ and $S(N, h)$.  Here $S(M, g)$ and $S(N, h)$ are the 
unit tangent bundles of $(M, g)$ and $(N, h)$ respectively.  We call $F$ a 
$C^{k}$-{\it geodesic conjugacy} from $M$ to $N$. A compact Riemannian manifold
is
said to be $C^{k}$-{\it geodesically rigid}  within a given class of manifolds if any Riemannian manifold
$N$ in that class whose geodesic flow 
is $C^{k}$-conjugate to that of $M$ is isometric to $M$.

A. Weinstein (\cite{W}) exhibited a zoll surface of non-constant curvature
whose
geodesic flow is  conjugate to that of the round sphere.  On 
the other hand, two flat tori with $C^{1}$-conjugate geodesic flows must be 
isometric. Therefore, a natural question arises:

\proclaim{Question}
Which compact Riemannian manifolds are $C^{k}$-geodesically rigid in a given
class of manifolds?
\endproclaim

This question is central to the study of negatively curved manifolds.
Many important 
open problems in the field will follow if all negatively curved manifolds can
be shown to be
geodesically rigid (see \cite{BFL}, \cite{EHS}, \cite{Ka},
\cite{Kk}). For negatively
curved surfaces , C. Croke (\cite{C}) and J. Otal (\cite{O}) have 
independently answered this question affirmatively. Recently, C. Croke and B.
Kleiner (\cite{CK})
 proved that compact Riemannian manifolds with a parallel vector field are
geodesically rigid.
 For negatively curved manifolds of higher dimension, the question is still
open. Recently, G. Besson,
G. Courtois and G. Gallot proved that if a manifold has a geodesic flow
which is $C^{1}$-conjugate to the geodesic flow
of a rank one locally symmetric  space $M$ , then
it is isometric to $M$.

In this article, we will consider the question above for the class of \rom{2}-step
Riemannian nilmanifolds.
This question was also asked by P. Eberlein in \cite{E}. 

A $k$-{\it step
Riemannian nilmanifold}
 is a compact quotient of a $k$-step nilpotent Lie group $G$ by a discrete
subgroup together with a 
Riemannian metric whose lift to $G$ is left-invariant.  Note that a one-step
nilmanifold is just a 
flat torus.  Thus the \rom{2}-step nilmanifolds may be viewed as the simplest
generalization of flat tori,
 yet they have a much richer geometry.

We are interested in this problem in part because of its relationship to
spectral geometry. 
 Two Riemannian manifolds are said to be {\it isospectral} if the associated
Laplace-Beltrami
 operators have the same eigenvalue spectrum.  A continuous family $M_{t}$ 
($-\epsilon 
< t < \epsilon$) of Riemannian manifolds is said to be an {\it isospectral
deformation} of $M_0$ 
if the manifolds are pairwise isospectral.  Since the Laplacian can be viewed
as the quantum analog of the 
geodesic flow, one might expect that continuous families of isospectral
manifolds would have conjugate geodesic flows.
E. N. Wilson and the first author (\cite{GW}) gave a method for constructing
isospectral
 deformations of Riemannian nilmanifolds, in particular, \rom{2}-step nilmanifolds. 
 Moreover, P. Eberlein (\cite{E}) showed that if a pair of \rom{2}-step
nilmanifolds have $C^1$ conjugate
 geodesic flows, then they must both lie in one of these continuous families. 
Nonetheless, we shall see
that, at least for some of the isospectral deformations, the geodesic flows are
not even $C^{0}$ conjugate.  

The main results in this paper are the following:

\proclaim{Theorem 1}

Almost all compact \rom{2}-step nilmanifolds are $C^{2}$-geodesically rigid within
the class of all compact nilmanifolds.
\endproclaim

(See section three for a precise statement.)

\smallskip

\proclaim{Theorem 2}
There is a large class of compact \rom{2}-step Riemannian nilmanifolds such that
any Riemannian
 nilmanifold whose geodesic flow is $C^0$ conjugate to a nilmanifold $M$ in
this class is isometric to $M$.
Included in this class are many continuous families of isospectral manifolds.
\endproclaim

\smallskip

We remark that R. Kuwabara \cite{Ku} has shown that for some of the
isospectral deformations
 of \rom{2}-step nilmanifolds mentioned above, the geodesic flows restricted to
dense open subsets
 of the tangent bundles are symplectically conjugate.

A. Kaplan \cite {Kap} introduced the notion of \rom{2}-step nilmanifolds of
Heisenberg type.  (See section one for the definition.)  

\proclaim{Theorem 3}
Compact \rom{2}-step nilmanifolds of Heisenberg type are $C^{0}$-geodesically rigid
within the class of compact 
Riemannian nilmanifolds.
\endproclaim

\smallskip

Finally, Eberlein \cite{E} defined a notion of Riemannian \rom{2}-step
nilmanifolds in resonance. 
 This algebraic condition is closely associated
with the geometric condition that the set of vectors lying in closed orbits of
the geodesic
 flow are dense in the unit tangent bundle.  (See \cite {Ma}, \cite{LP}.)  We define a
notion of 
nilmanifold {\it  strongly in resonance} and prove:

\proclaim{Theorem 4}
Compact \rom{2}-step nilmanifolds which are strongly in resonance are
$C^{2}$-geodesically rigid within the class of compact 
Riemannian nilmanifolds.
\endproclaim

\smallskip

The paper is organized as follows:  In section one, we review the geometry of
\rom{2}-step nilmanifolds, 
define various special classes of \rom{2}-step nilmanifolds, including those cited
in the theorems above, and 
review the construction of isospectral deformations of nilmanifolds.  Section
two establishes a number of
results concerning conjugacies
of geodesic flows between arbitrary \rom{2}-step nilmanifolds and culminates in
Theorem 2.  Theorems 1, 3 and 4
are established in section three.

The authors would like to thank Patrick Eberlein, Jens Heber, and David Webb
for helpful discussions.  Some of the 
results of this paper were announced in \cite {GM}.

\heading 1. Preliminaries\endheading

In this section, we give a brief introduction to the geometry of \rom{2}-step Lie
groups that will be used in the proof of the theorems. We recommend
\cite{E} as a reference.

A Riemannian nilmanifold is a quotient $M = \Gamma\backslash N$ of a nilpotent
Lie group $N$ by a discrete subgroup $\Gamma$, together with a Riemannian
metric $g$ whose lift to $N$, also denoted $g$, is left-invariant.
We say the nilmanifold has step size $k$ if $N$ is $k$-step nilpotent.

\proclaim{1.1 Notation}
Let $N$ be a \rom{2}-step nilpotent Lie group with a left invariant metric $g$.
The metric $g$ defines an inner product $< , >$ on the Lie algebra $\Cal N$
of $N$. Let $\Cal{Z} = [\Cal{N}, \Cal{N}]$ and let $\Cal V$ denote the orthogonal complement of $\Cal Z$ in 
$\Cal N$ relative to $< , >$.

Since $N$ is \rom{2}-step nilpotent, the Campbell-Baker-Hausdorff theorem gives the product rule
$$
\exp(x)\exp(y) = \exp(x + y + \frac{1}{2}[x, y])
$$ 
for all $x, y \in {\Cal N},$
where $\exp: {\Cal N} \rightarrow N$ is the Lie group exponential map.

For $z$ in $\Cal Z$, define a skew symmetric linear transformation  
$J(z): {\Cal V} \rightarrow {\Cal V}$  by 
$J(z)x = (\operatorname{ad}(x))^{*}z$ for $x\in \Cal{V}$,
where $(\operatorname{ad}(x))^{*}$ denotes the adjoint of $\operatorname{ad}(x)$. Or, equivalently

\roster
\item"($\ast$)" $\qquad\qquad <J(z)x, y> = <[x, y], z> \quad\text{for $x, y \in {\Cal V}, z \in {\Cal Z}$.}$
\endroster

\noindent The operators $J(z)$  contain all geometric information concerning the manifold $(N, g).$

Conversely, given inner product spaces $\Cal{V}$ and $\Cal{Z}$ and a linear map
$J: \Cal{Z}\rightarrow {\frak so}(\Cal{V})$, we can construct a \rom{2}-step
nilpotent Lie algebra $\Cal{N}$ together with an inner product, by setting
$\Cal{N} = \Cal{V}\bigoplus\Cal{Z}$ and defining the Lie bracket so that $\Cal{Z}$ is central and $[\quad,\quad]: \Cal{V}\times\Cal{V}\rightarrow
\Cal{Z}$ is given by $(\ast)$. The inner product defines a left-invariant Riemannian metric on the associated simply-connected nilpotent Lie group $N$.
\endproclaim

\smallskip

\proclaim{1.2 Example}
The Heisenberg group of dimension $2n + 1$ is the simply-connected Lie group
with Lie algebra $\Cal{N} = \operatorname{span}\{ x_1, \cdots, x_n, y_1, 
\cdots, y_n, z\}$, where $[x_i, y_i] = z, \quad 1\leq i \leq n, $ and all
other brackets of basis elements are zero. Giving $\Cal{N}$ the inner product 
for which the basis above is orthonormal, we have 
$$
\align
\Cal{Z} &=
\operatorname{span}\{z\},\\
\Cal{V} &= \operatorname{span}\{ x_1, \cdots, x_n, y_1, \cdots, y_n\}
\endalign
$$
and $J(z)x_i = y_i,\quad J(z)y_i = - x_i.$ Thus $J(z)$
can be viewed as a complex structure on the vector space $\Cal{V}$.
\endproclaim

\proclaim{1.3 Definition}
\roster
\item A \rom{2}-step nilpotent Lie group $N$ with a left-invariant Riemannian metric $g$ is of {\it Heisenberg }type if
$ 
J^{2}(z) = -|z|^{2} \Id  \quad\text {on $\Cal V$, for all $z\in {\Cal Z}$}.
$

\item $(N, g)$ is said to be {\it nonsingular} if $J(z)$ is a nonsingular linear transformation
for all $z\in {\Cal Z}$, $z \neq 0$.

\item $(N, g)$ is said to be {\it in resonance} if for 
all $z\neq 0 
\in {\Cal Z}$,  every ratio of non-zero eigenvalues of $J(z)$ is rational, 
equivalently, for each $z\in {\Cal Z}$, there exists a constant $t$ such 
that $e^{tJ(z)} = \operatorname{\Id }$, where $t$ may depend on $z$.

\item $(N, g)$ is said to be {\it  strongly in resonance} if for each $z \in 
{\Cal Z}$, there exists a constant  $t$ such that $e^{tJ(z)} = 
-\operatorname{\Id }$, where $t$ may depend on $z$.

\item $(N, g)$ is said to be {\it irrational} if for $z$ in a dense set of 
${\Cal Z}$, $1,\, \theta_{1}(z),\, \cdots,\, \theta_{p}(z)$ are linearly independent 
over the field of rational numbers, where $\{\pm \sqrt{-1}\theta_{i}(z),\, 
i= 1, 2, \cdots, p \}$ are the distinct eigenvalues of $J(z)$. \rom{(} We choose $\theta_{i}(z)
\geq 0$\rom{)}.
\endroster
\endproclaim

\proclaim{1.4 Remark} 
\roster
\item  The notion of Lie group of Heisenberg type is due to A. Kaplan
(\cite{Kn}). These groups have a surprisingly rich and varied geometry and 
have been studied by several authors (e.g., \cite{E}, \cite{Rie}). It is easy to see
that the Lie groups of Heisenberg type are nonsingular and strongly in resonance.

\item P. Eberlein (\cite{E}) defined the term ``in resonance''. 
The nilmanifolds of Heisenberg type are in resonance, in fact, strongly in resonance. M. Mast (\cite{Mt}) gave the first examples of nilmanifolds in resonance which are not of Heisenberg type. We give an additional example below. 
\item The condition ``strongly in resonance'' implies ``in resonance''. 
However, our 
examples  will show that ``strongly in resonance'' is indeed stronger than 
``in resonance''.

\item Generically, all \rom{2}-step nilpotent Lie groups are irrational. In this 
case, for $z$ in a dense set of ${\Cal Z}$, there is a sequence of $t_{i}$
such that $\lim_{i\rightarrow \infty}e^{t_{i}J(z)} = -\Id .$ 
\endroster
\endproclaim

\smallskip

\proclaim{1.5 Example} Let $\Cal{V}$ and $\Cal{Z}$ be inner product spaces with orthonormal bases $\{X_1,\cdots, X_8\}$ and $\{Z_1, Z_2\}$, respectively, and define $J: \Cal{Z}\rightarrow {\frak so}(\Cal{V})$ so that for 
$Z = z_1 Z_1 +  z_2 Z_2$, $J(Z)$ has the following matrix representation with respect to the basis above of $\Cal{V}$:
$$
\left(
\matrix
0&0&-\lambda_1 z_{1}&-\lambda_1 z_{2}&0&0&0&0& \\
0&0&\lambda_{1} z_{2}&-\lambda_{1} z_{1}&0&0&0&0 \\
\lambda_{1} z_{1}&-\lambda_{1} z_{2}&0&0&0&0&0&0 \\
\lambda_{1} z_{2}&\lambda_{1} z_{1}&0&0&0&0&0&0 \\
0&0&0&0&0&0&-\lambda_2 z_{1}&-\lambda_2 z_{2} \\
0&0&0&0&0&0&\lambda_{2} z_{2}&-\lambda_{2} z_{1}\\
0&0&0&0&\lambda_{2} z_{1}&-\lambda_{2} z_{2}&0&0\\
0&0&0&0&\lambda_{2} z_{2}&\lambda_{2} z_{1}&0&0 \\
\endmatrix
\right)
$$
The distinct eigenvalues of $J(Z)$ are 
$\{ \pm \sqrt{-1} \lambda_{1}|Z|,\, \pm \sqrt{-1} \lambda_{2}|Z|\}$.

As in \rom{1.1}, the data $(\Cal{V}, \Cal{Z}, J)$ defines a simply-connected \rom{2}-step nilpotent Lie group $N$ with a left-invariant Riemannian metric $g$. If, say, 
$\lambda_{2} = 2\lambda_{1}$, then $(N, g)$ is in resonance but not strongly in resonance. If, say, $\lambda_{2} = 3\lambda_{1}$, then $(N, g)$ is strongly in
resonance but is not of Heisenberg type.
\endproclaim

We now consider geodesics on a 2-step nilmanifold. 
Let $x, y \in {\Cal N}$ regarded as left invariant vector fields on $N$. 
Recall that the formula for the covariant derivative $\bigtriangledown_{x}y$
normally contains 6 terms but in this case three of them vanish since $g(x, y)$
is constant in $N$. One obtains
$$
\bigtriangledown_{x}y = \frac{1}{2}\{[x, y] -  (\operatorname{ad}(x))^{*}y 
- (\operatorname{ad}(y))^{*}x\}
$$
Therefore, one obtains the following lemma.

\proclaim{1.6 Lemma}\rom{(\cite{E})}

\roster
\item  $\bigtriangledown_{x}y = \frac{1}{2}[x, y]$, \quad for  $x, y \in 
{\Cal V}$.

\item  $\bigtriangledown_{x}z = \bigtriangledown_{z}x = -\frac{1}{2}J(z)x,$
       \quad for  $x \in {\Cal V}, z \in {\Cal Z}$.

\item  $\bigtriangledown_{z}z^{*} = 0$, \quad for  $z, z^{*} \in {\Cal Z}.$
\endroster
\endproclaim

\smallskip

Using lemma 1.6,  one can calculate geodesics in $N$.

Let $\sigma(t)$ be a geodesic with $\sigma(0) = e$ and $\sigma'(0) =
x + z \in {\Cal N}$, where $e$ is the identity element of $N$, 
$x \in {\Cal V}$ and $z \in {\Cal Z}$. In exponential coordinates, we write
$\sigma(t, x + z) = \exp(x(t) + z(t))$ where $x(t) \in {\Cal V}, z(t) 
\in {\Cal Z}, x'(0) = x, z'(0) = z.$  We denote $J(z)$ by $J$.  One obtains:

\smallskip

\proclaim{1.7 Lemma}\rom{(\cite{E})}
\roster
\item  $\sigma'(t) = dL_{\sigma(t)}(e^{tJ}x + z)$,\quad for  $t 
   \in {\Bbb R},$  

\noindent where $L_{m}$ denotes the left translation on $N$ by $m \in N$.

\item $x(t) = tx_{1} + (e^{tJ} - \id ) J^{-1}x_{2}$,\quad for  $t\in {\Bbb R}$,

   \noindent where $x_{1} \in {Ker}(J(z)), x_{2}\in {Ker}(J(z))^{\perp}$
    and $x_{1} + x_{2} = x.$

\item 
$z(t) = tz_{1}(t) + z_{2}(t), $
$$
\align
z_{1}(t) & = z + \frac{1}{2}[x_{1}, (e^{tJ} + \Id )J^{-1}x_{2}] + 
                 \frac{1}{2}\sum_{i=1}^{N}[J^{-1}\xi_{i}, \xi_{i}],  \\
\split z_{2}(t) & =[x_{1}, (\Id  - e^{tJ})J^{-2}x_{2}] -  
                    \frac{1}{2}[e^{tJ}J^{-1}x_{2}, J^{-1}x_{2}]   \\
         &-\frac{1}{2}\sum_{i\neq j=1}^{N} \frac{1}{\theta_{j}^{2}-
       \theta_{i}^{2}}\{[e^{tJ}J\xi_{i}, e^{tJ}J^{-1}\xi_{j}]     
         - [e^{tJ}J\xi_{i}, e^{tJ}\xi_{j}]\}   \\
   &+\frac{1}{2}\sum_{i\neq j=1}^{N}\frac{1}{\theta_{j}^{2}-
     \theta_{i}^{2}}
     \{[J\xi_{i}, J^{-1}\xi_{j}] - [\xi_{i}, \xi_{j}]\},\endsplit
\endalign
$$
where $e^{tJ} = \sum_{n}\frac{t^{n}}{n!}J^{n}$, 
$\{\pm \theta_{i}\sqrt{-1}, ; i=1, \cdot\cdot\cdot, p \}$ are the distinct 
eigenvalues 
of $J(z)$ and $\{\xi_{j}\} \subseteq{Ker}(J(z))^{\perp}$ satisfies 
$\Sigma_{j} \xi_{j} = x_{2}$ and $J^{2}(z) \xi_{j} = - \theta^{2}_{j} \xi_{j}$. 
\endroster
\endproclaim

\proclaim{1.8 Remark}
\roster
\item
If $N$ is of Heisenberg type, then $x_{1} = 0,$ and $\pm \sqrt{-1}|z|$ 
are the only eigenvalues of $J(z)$, so
$$\align
        x(t)  &= (\cos(t\mid z\mid) - 1)J^{-1}x + \frac{\sin(t\mid z\mid)
                   }{\mid z\mid} x,   \\
       z_{1}(t) &= (1 + \frac{\mid x\mid^{2}}{2\mid z\mid^{2}})z,  \\
       z_{2}(t) &= -\frac{\sin(t\mid z\mid) \mid x\mid^{2}}{2\mid z\mid^{3}}z 
\endalign
$$
\item
If $x = 0$ or $z = 0$, then the geodesics become respectively
$\sigma(t) = \exp(tz)$ or $\sigma(t) = \exp(tx).$ More generally, if
$J(z)x = 0,$ then the geodesic is given by $\sigma(t) =
\exp(t(x + z)).$
\endroster
\endproclaim

\smallskip

A geodesic $\sigma$ in $N$ descends to a closed geodesic in $\Gamma\backslash N$
if and only if there exists $\gamma\in \Gamma$ and $l\in \Bbb{R}$ such that 
$\sigma(t + l) = \gamma \sigma(t)$ for all $t\in\Bbb{R}.$

\proclaim{1.9 Lemma}\rom{(\cite{E})}
Let $\varphi = \exp(v^{*} + z^{*})$ be an arbitrary element of $N$. Let
$z^{**}$ be the component of $z^{*}$ orthogonal to $[v^{*},\, {\Cal V}].$
Let $\sigma(t)$ be a unit speed geodesic such that $\varphi \sigma(t)
= \sigma(t + l)$ for all $t \in R$ and some $l > 0.$
Let $a = \sigma(0)$ and $\omega^{*} = ( |v^{*}|^{2} + 
|z^{**}|^{2})^{\frac{1}{2}}.$ Then
\roster
\item $|v^{*}| \leq l \leq \omega^{*}.$

\item $l = \omega^{*}$ if and only if the following conditions hold:
$$
\align
\sigma(t) &= a \exp(\frac{t}{\omega^{*}}(v^{*} 
                     + z^{**}) \\
z^{**} &= z^{*} + [v^{*}, \xi],\quad \text{ where} \,\,a = \exp(\xi).
\endalign
$$
\item $l = |v^{*}|$ if and only if $z^{**} = 0$ and  $\sigma(t) 
         = a \cdot \exp(t \frac{v^{*}}{|v^{*}|}).$
\endroster
\endproclaim

\smallskip

We next recall some basic facts about uniform discrete subgroups of $N$. See
\cite{Rn} for details.

\proclaim{1.10} If $\Gamma$ is a uniform discrete subgroup of $N$, then $\log\Gamma\cap\Cal{Z}$ is a lattice of full rank in the derived algebra $\Cal{Z}$, and $\pi_{v}(\log\Gamma)$ is a lattice of full rank in $\Cal{V}$, where $\pi_{v}: \Cal{N}\rightarrow \Cal{V}$ is the projection with kernel $\Cal{Z}$. For $x, y \in \log\Gamma$, we have 
$\exp(x)\exp(y)\exp(-x)\exp(-y)=\exp([x, y])\in \Gamma$, so $[x, y] \in \log\Gamma\cap\Cal{Z}$. In particular, if we choose a basis of $\Cal{N}$ consisting of elements of $\log\Gamma$, then the constants of structure are rational. Thus, letting $\Cal{N}_{Q}$ denote the rational span of $\log\Gamma$, then $\Cal{N}_{Q}$ has the structure of a rational Lie algebra. We will say a subspace of $\Cal{N}$ is {\it rational} if it has a basis consisting of elements of $\log\Gamma$. For example, $\operatorname{image}(\operatorname{ad}(x))$ is rational for all $x\in \log\Gamma$.
\endproclaim

We conclude this section with a brief discussion of isospectral nilmanifolds. The notion of almost inner automorphisms (1.11 below) will play a key role in the sequel.

\proclaim{1.11 Definition}
\roster
\item"(a)"  Let $\Gamma$  be a uniform discrete subgroup of a 
simply-connected
nilpotent Lie group $N$.  An automorphism $\Phi$  of $N$ is said to be   
$\Gamma$-{\rm almost inner} if  $\Phi(\gamma)$ 
is conjugate to $\gamma$  for all $\gamma \in \Gamma$. The automorphism is 
said to be {\rm almost inner} if $\Phi(x)$ is conjugate to $x$ for all $x 
\in N$ .
\item"(b)" A derivation $\phi$ of the Lie algebra $\Cal{N}$ is said to be $\Gamma$-almost inner, respectively almost inner, if $\phi(x) \in \operatorname{image}(\operatorname{ad}(x))$ for all $x\in \log\Gamma$, respectively,  for all $x\in \Cal{N}$.
\endroster
\endproclaim

\proclaim{1.12 Remark} \rom{(See \cite{GW})}
\roster
\item"(a)" The $\Gamma$-almost inner automorphisms and the almost inner automorphisms form connected Lie subgroups of $\operatorname{Aut}(N)$. In many cases, these groups properly contain the group $\operatorname{Inn}(N)$ of inner automorphisms. The spaces of $\Gamma$-almost inner \rom{(}respectively, almost inner\rom{)} derivations of $\Cal{N}$ are the Lie algebras of these groups of automorphisms. In particular, if $\phi$ is a \rom{(}$\Gamma$-\rom{)}almost inner derivation, then there exists a one-parameter family $\Phi_{t}$ of \rom{(}$\Gamma$-\rom{)}almost inner automorphisms of $N$ such that ${\Phi_{t}}_{*} = e^{t\phi}$. Conversely, if $\Phi$ is a \rom{(}$\Gamma$-\rom{)}almost inner automorphism of $N$, then $\Phi_* = e^{\phi}$ for some 
\rom{(}$\Gamma$-\rom{)}almost inner derivation of $\Cal{N}$.
\item"(b)" Note that a $\Gamma$-almost inner derivation $\phi$ satisfies $\phi(\Cal{V}) \subseteq [\Cal{N}, \Cal{N}]$ and $\phi(Z) = 0$ if $Z$ is central. In particular, if $\Cal{N}$ is \rom{2}-step nilpotent, then (letting 
$\Cal{Z} = [\Cal{N}, \Cal{N}]$ as before), we have $\phi(\Cal{N}) \subseteq \Cal{Z}$ and $\phi(\Cal{Z}) = 0$, so $\phi^{2} = 0.$ Thus $e^{t\phi} =
\operatorname{\Id } + t\phi.$
\endroster
\endproclaim

\proclaim{1.13 Proposition}
Let $(\Gamma\backslash N, g)$ be a compact nilmanifold and let $\Phi$ be a $\Gamma$- almost inner automorphism of $N$. Then $(\Phi(\Gamma)\backslash N, g)$
is isospectral to $(\Gamma\backslash N, g)$. Conversely, if $N$ is \rom{2}- step nilpotent and if $(\Gamma_{t})_{t\geq 0}$ is a continuous family of discrete subgroups of $N$ such that the family of manifolds $(\Gamma_t\backslash N, g)$
are all isospectral, then there exists a family $\{\Phi_t\}$ of $\Gamma_0$-almost inner automorphisms of $N$ such that $\Gamma_t = \Phi_t(\Gamma_0)$ for all $t$.
\endproclaim

The first statement is proven in \cite{GW} for almost inner automorphisms and in \cite{G} for $\Gamma$-almost inner automorphisms. The converse is given in \cite{OP} and \cite{P}.

\proclaim{1.14 Remark}
If $\Phi$ is an inner automorphism of $N$, say $\Phi$ is conjugation by 
$a\in N$, then $(\Phi(\Gamma)\backslash N, g)$ is isometric to 
$(\Gamma\backslash N, g)$. The isometry is induced from the isometry $L_a$ of
$(N, g)$ given by left translation. However, if $\phi$ is a $\Gamma$-almost inner derivation which is not inner and $\{\Phi_t\}$ is the corresponding family of automorphisms, then the deformation $\Phi_{t}(\Gamma)\backslash N, g)$ is non-trivial.
\endproclaim

\heading 2. Geodesic Conjugacy \endheading

Throughout this section, $(N, g)$ will denote an arbitrary simply-connected 2-step nilpotent Lie group with left-invariant Riemannian metric and $\Gamma$ a uniform discrete subgroup of $N$. We continue to use the notation $\Cal{N} =
\Cal{Z} + \Cal{V}$ introduced in 1.1 for the Lie algebra of $N$.

\proclaim{2.1 Notation and Remarks}
\roster
\item"(a)"The left-invariant vector fields on $N$ induce global vector fields on $\Gamma\backslash N$. Thus the tangent bundles of both $N$ and 
$\Gamma\backslash N$ are completely parallelizable, and the unit tangent bundles may be identified with 
$$
S(N, g) = N \times S(\Cal{N})
$$
$$
S(\Gamma\backslash N, g) = \Gamma\backslash N \times S(\Cal{N})
$$
where $S(\Cal{N})$ is the unit sphere in $\Cal{N}$ relative to the Riemannian inner product. We will write $S(N)$ for $S(N, g)$ if $g$ is understood.

\item"(b)" For $x\in N$, the left action $L_x: N\rightarrow N$ induces a diffeomorphism $(L_x)_*$ of $S(N)$. Under the identification above, 
$$
(L_x)_*(n, u) = (xn, u)
$$ 
for $n\in N, u\in S(\Cal{N}).$

\item"(c)" Now suppose that $(\Gamma^*\backslash N^*, g^*)$ is another 
\rom{2}-step nilmanifold and that
$$
F: S(\Gamma\backslash N, g) \rightarrow S(\Gamma^*\backslash N^*, g^*)
$$
is a homeomorphism intertwining the geodesic flows. By \rom{(a)}, the fundamental group of $ S(\Gamma\backslash N)$ is isomorphic to $\Gamma$, so $F$ induces an isomorphism 
$$
F_{*}: \Gamma \rightarrow \Gamma^*.
$$
\endroster
\endproclaim

\smallskip

\proclaim{2.2 Proposition} \rom{(\cite{E})}
Suppose $(\Gamma\backslash N, g)$ and $(\Gamma^*\backslash N^*, g^*)$ are compact \rom{2}-step nilmanifolds and $F: S(\Gamma\backslash N)\rightarrow
S(\Gamma^*\backslash N^*)$ is a homeomorphism intertwining their geodesic flows. Then there exists a $\Gamma$-almost inner automorphism $\Phi$ of $N$ (see definition 1.11) such that $(\Gamma^*\backslash N^*, g^*)$ is isometric to 
$(\Gamma\backslash N, g)$. Moreover the isomorphism $F_{*}: \Gamma \rightarrow
\Gamma^*$ is given by $F_* = \Psi_* \circ \Phi|_{\Gamma}$ where
$$
\Psi: (\Psi(\Gamma)\backslash N, g) \rightarrow (\Gamma^*\backslash N^*, g^*)
$$
is an isometry and $\Psi_{*}: \Phi(\Gamma) \rightarrow \Gamma^*$
is the induced map on fundamental groups.
\endproclaim

\smallskip

\proclaim{2.3 Notation and Remarks}
\roster
\item"(a)" By \rom{Proposition 2.2}, we may replace $(\Gamma^*\backslash N^*, g^*)$ 
by $(\Phi(\Gamma)\backslash N, g)$, so $F: S(\Gamma\backslash N, g)\rightarrow
S(\Phi(\Gamma)\backslash N, g)$, and we may assume that $F_* = \Phi|_{\Gamma}$ for some $\Gamma$-almost inner automorphism $\Phi$ of $N$. Let 
$$
\Tilde{F}: S(N) \rightarrow S(N)
$$
be the lift of $F$. In the notation of \rom{2.1 (b)}, we have 
$$
\Tilde{F}\circ {L_{\gamma}}_* =  {L_{\Phi(\gamma)}}_* \circ \Tilde{F},
$$
for all $\gamma \in \Gamma.$
Also, denoting by $G^t$ the geodesic flow of $(N, g)$, we have
$$
\Tilde{F}\circ G^t = G^t \circ \Tilde{F}.
$$

\item"(b)" For $(n, u)\in S(N)$, write 
$$
\Tilde{F}(n, u) = (\exp(A(n, u) + B(n, u))n, \,\,
 I(n, u) + H(n, u))
$$

\item"(c)" \rom{(}See \rom{Remark 1.12}\rom{)} There exists a $\Gamma$-almost inner derivation of $\Cal{N}$, which we denote by $\phi$, such that the differential 
$\Phi_{*}: \Cal{N} \rightarrow \Cal{N}$ is given by $\Phi_{*} = \operatorname{\Id } + \phi$.  We have $\phi({\Cal Z}) = 0$ and $\phi(\Cal{V})\subseteq \Cal{Z}.$
\endroster
\endproclaim

\smallskip

\proclaim{2.4 Proposition} We use the notation in \rom{2.3}.
For $(n, u) \in S(N)$ and $\gamma\in\Gamma$,
\roster
\item"(a)"  $B(\gamma n, u) = B(n, u),$

\item"(b)"  $H(\gamma n, u) = H(n, u),$

\item"(c)"  $I(\gamma n, u)  = I(n, u),$

\item"(d)"  $A(\gamma n, u) = A(n, u) + \phi(\log\gamma) 
     -[B(n, u), \log\gamma].$
\endroster
\endproclaim

\smallskip

\demo{Proof}  By \rom{2.1 (b)} and \rom{2.3}, we have 
$$
\align
&(\exp(A(\gamma n, u) + B(\gamma n, u))\gamma n, \,\,
 I(\gamma n, u) + H(\gamma n, u))  \\
&=\Tilde{F}(\gamma n, u) = {L_{\Phi(\gamma)}}_* \Tilde{F}(n, u)  \\
 &=(\Phi(\gamma)\exp(A(n, u) + B(n, u))n, \,\,
 I(n, u) + H(n, u)).
\endalign
$$
This yields \rom{(b)} and \rom{(c)}. Using \rom{2.3 (c)} and the fact that
$$
\exp(x)\exp(y) = \exp(x + y + \frac{1}{2} [x, y])
$$ for $x, y \in \Cal{N},$ a straight-forward computation yields \rom{(a)} and \rom{(d)}.
\enddemo

\smallskip

\proclaim{2.5 Proposition}
For a unit vector $u\in\Cal{N}$, let $\sigma(t, u)$ be the geodesic in $N$ with $\sigma(0, u) = e$ and $\dot{\sigma}(0, u) = u$. Write 
$\sigma(t, u) = \exp(X(t, u) + Z(t, u))$ with 
$X(t, u) \in \Cal{V}$ and $Z(t, u) \in \Cal{Z}$. Defining $J: \Cal{Z}\rightarrow {\frak so}(\Cal{V})$ as in \rom{1.1} and letting $G^t$ be the geodesic flow of $(N, g)$, we have:
\roster
\item"(a)" $I(G^{t}(n, u))  = e^{tJ(H(n, u))} I(n, u), $

\item"(b)" $H(G^{t}(n, u))  = H(n, u),  $

\item"(c)" $B(G^{t}(n, u)) = B(n, u) + X(t, I(n, u)+H(n, u))
                      - X(t, u) $

%\item"(d)" $A(G^{t}(n, u)) + \frac{1}{2}[B(G^{t}(n, u)), \log(n)]$ 
                   
% $ + \frac{1}{2}[B(G^{t}(n, u)), X(t, u)] $
%$  + \frac{1}{2}[\log(n), X(t, u)] + Z(t, u) $

%\noindent $= A(n, u)$
%$+ \frac{1}{2}[B(n, u), \log(n)] 
%+ \frac{1}{2}[B(n, u), X(t, I(n, u) + H(n, u))] $ 

%$ + \frac{1}{2}[\log(n), X(t, I(n, u) + H(n, u))] 
%                  + Z(t, I(n, u) + H(n, u)). $
\endroster
\endproclaim

\smallskip

\demo{Proof} The proposition follows from the expression for the geodesics in \rom{Lemma 1.7}.
\enddemo

As seen in Lemma 1.7, if $u\in \Cal{V}$ or $u\in \Cal{Z}$, then the geodesic $\sigma(t, u)$ in $N$ with $\sigma(0, u) = e, \dot{\sigma}(0, u) = u$ is given by $\sigma(t, u) = \exp(tu).$ For $n\in N,$ the geodesic with initial conditions $(n, u)$ is just the translation $n\exp(tu).$ The corresponding orbit of the geodesic flow in $N$ is the curve $(n\exp(tu), u)$. 

\proclaim{2.6 Proposition} In the notation of \rom{1.1} and \rom{2.3}, we have:
\roster
\item"(a)" $F: \Gamma\backslash N \times \{z\}\rightarrow \Phi(\Gamma)\backslash N \times \{z\} $ for all $z\in \Cal{Z}$.

\item" " That is, $H(n, z) = z$ and $ I(n, z) = 0$.

\item"(b)" $F: \Gamma\backslash N \times \{v\}\rightarrow \Phi(\Gamma)\backslash N \times \{v\} $ for all $v\in \Cal{Z}$.

\item" " That is, $H(n, v) = 0$ and $ I(n, v) = v$.
\endroster
\endproclaim

\smallskip

\proclaim{2.7 Lemma} \rom{(\cite{E})} Let $v\in\pi_{v}(\log\Gamma)$ and
 $\epsilon > 0$ be given. Then there exists an element $\xi \in \log\Gamma$ such that $\xi = kv + z_{0}$ for some positive integer $k$ and some element 
$z_{0}\in {\Cal Z}$ with $|z_{0}| < \epsilon.$
\endproclaim
 
\demo{Proof of Proposition 2.6}
Recall that free homotopy classes of closed curves in $\Gamma\backslash N$ correspond to conjugacy classes in $\Gamma.$ For $\gamma \in \Gamma$ (respectively, $\Phi(\Gamma$)), we will denote by $[\gamma]$ (respectively, $[\gamma]^*$) the corresponding free homotopy class.

\smallskip

\noindent (a) Let
$$
\align
\Cal{Z}_{u} &= \{z \in \Cal{Z}, \|z\| = 1\} \\ 
\Cal{Z}_{\Gamma} &=
\{z \in \Cal{Z}_{u}\, | rz \in \log(\Gamma),\, \text{for some}\, r > 0\}
\endalign
$$ 
Since $\Gamma$ is a uniform lattice, $\Cal{Z}_{\Gamma}$ is dense
in $\Cal{Z}_{u}$. Thus  we  need only  show that $H(n, z) = z$ for $z \in \Cal{Z}_{\Gamma}.$

The almost inner automorphism $\Phi$ restricts to the identity on the center of $\Cal{N}$. In particular, $\log\Gamma\cap\Cal{Z} = \log\Phi(\Gamma)\cap\Cal{Z}$ and $\Cal{Z}_{\Gamma} = \Cal{Z}_{\Phi(\Gamma)}.$ Let $z_{0} \in \Cal{Z}_{\Gamma}$
and let $r_{0}$ satisfy $r_0 z_0 \in \log\Gamma$. By \rom{Lemma 1.9} (with $v^* = 0$), the longest geodesics in the free homotopy class $[r_0 z_0]$ are precisely the projections to $\Gamma\backslash N$ of $\exp(tz_0)$ and all its left translations. Thus the submanifold $\Gamma\backslash N \times\{z_0\}$ of $S(\Gamma\backslash N)$ is foliated by all the longest periodic orbits of the
geodesic flow in the free homotopy class $[r_0 z_0]$ (viewed now as a free homotopy class of curves in  $S(\Gamma\backslash N)$). Similarly 
$\Phi(\Gamma)\backslash N \times\{z_0\}$ is foliated by the longest periodic orbits of the flow in the class $[r_0 z_0 ]^*$ of $S(\Phi(\Gamma)\backslash N)$.
Hence the geodesic conjugacy $F$ must map $\Gamma\backslash N \times\{z_0\}$
onto $\Gamma^*\backslash N \times\{z_0\}$. This proves \rom{(a)}.

\smallskip

\noindent (b) The proof of \rom{(b)} is similar but considerably more complicated.   Let 
$$
\align
\Cal{V}_{u} &= \{v \in \Cal{V} | \,\,\|v\| = 1\}   \\
\Cal{V}_{\Gamma} &= \{v \in \Cal{V}_{u} | rv \in \pi_{v}(\log\Gamma),
      \text{for some} \,\,r > 0\}
\endalign
$$
where $\pi_{v}: \Cal{N}\longrightarrow \Cal{V}$ is the orthogonal projection.
Since $\Gamma$ is a uniform lattice, $\Cal{V}_{\Gamma}$ is dense in
$\Cal{V}_{u}$ and  we only need to show that $I(n, v) = v$ for all 
$v\in{\Cal V}_{\Gamma}$. 

For a fixed $v \in \Cal{V}_{\Gamma}$, there is an $r > 0$,
such that $rv \in \pi_{v}(\log\Gamma).$
According to  Lemma 2.7, for $\epsilon > 0$, there is a positive integer
$k_{\epsilon}$ and an element $z_{\epsilon} \in z$ with $\|z_{\epsilon}\|
< \epsilon$ such that
$$
k_{\epsilon}rv + z_{\epsilon} \in \log(\Gamma).
\tag 2.6.1
$$
Decompose 
$$
z_{\epsilon} = z_{\epsilon, 0} + z_{\epsilon, 1}
\tag 2.6.2
$$
with
$z_{\epsilon, 0} \in \operatorname{Im}(\operatorname{ad}(v))^{\perp}$ and 
$z_{\epsilon, 1} 
\in \operatorname{Im}(\operatorname{ad}(v))$.
Let 
$$
v_{0} = \frac{k_{\epsilon}rv}{\|k_{\epsilon}rv + z_{\epsilon, 0}\|}
\tag 2.6.3
$$
$$
z_{0} = \frac{z_{\epsilon, 0}}{\|k_{\epsilon}rv + z_{\epsilon, 0}\|}
\tag 2.6.4
$$
For $n \in N$, the curve $n\exp t(v_0 + z_0)$ is a geodesic 
(see Lemma \rom{1.7}). Set
$$
   N_{0} = \{ n \in N \,|\, \gamma_{n}(t) = n\exp t(v_{0} + z_{0})
                  \, \text{descends to a closed geodesic in} \,
                   \Gamma\backslash N\}.
$$
We are going to prove that $N_{0}$ is dense in $N$.

For $n = \exp(x + w) \in N$, $ x \in \Cal{V}, w \in \Cal{Z}$,
we get,
$$
\align
\gamma_{n}(t) &= n\exp(v_{0} + z_{0}) = \exp(tv_{0} + x + w +\frac{t}{2}
[x, v_{0}] + tz_{0}),\\
\gamma_{n}(t + l) &= \exp((t+l)v_{0} + x + w +\frac{t+l}{2}[x, v_{0}] +
 (t +l)z_{0})\\
&=\exp(lv_{0} + l[x, v_{0}] + lz_{0})
   \exp(x + tv_{0} + w + \frac{t}{2}[x, v_{0}] + tz_{0})\\
&=\exp(lv_{0} + l[x, v_{0}] + lz_{0})\gamma_{n}(t).
\endalign
$$
Therefore, $\gamma_{n}(t)$ descends to a closed geodesic in $\Gamma\backslash
N$ if and only if 
$$
\exp(lv_{0} + l[x, v_{0}] + lz_{0}) \in \Gamma. \tag 2.6.5
$$
So, 
$$
N_{0} = \{ n = \exp(x + w) \in N \,|\, \exp(lv_{0} + l[x, v_{0}] + l z_{0}) 
\in \Gamma\,\, \text{for some}\, l\in \Bbb{R}\}.
$$
By \rom{(2.6.1) - (2.6.4)}, we have  
$$
\|k_{\epsilon}rv + z_{\epsilon, 0}\|v_{0} + \|k_{\epsilon}rv + 
z_{\epsilon, 0}\|z_{0} + z_{\epsilon, 1} \in \log\Gamma.
$$
Since $z_{\epsilon, 1} \in \operatorname{Im}(\operatorname{ad}(v)),$ there is a 
$u\in \Cal{V}$ such that
$[u, v_{0}] = \frac{z_{\epsilon, 1}}{\|k_{\epsilon}rv + z_{\epsilon, 0}\|}$.
Letting $r_{1} = \|k_{\epsilon}rv + z_{\epsilon, 0}\|$, we get
$$
r_{1}(v_{0} + [u, v_{0}] + z_{0}) \in \log(\Gamma)
\tag 2.6.6
$$
Let $\Cal{Z}_{Q}$ be the rational span of $\Cal{Z}\cap \log(\Gamma)$; i.e., 
$\Cal{Z}_{Q} = \Cal{Z}\cap\Cal{N}_{Q}$ in the notation of \rom{1.10}. By the remark in \rom{1.10}, $\Cal{Z}_{Q}\cap 
\operatorname{Im}(\operatorname{ad}(v_{0}))$ is dense in
$\operatorname{Im}(\operatorname{ad}(v_{0}))$ and thus so is 
 $\frac{1}{r_{1}}(\Cal{Z}_{Q}\cap 
\operatorname{Im}(\operatorname{ad}(v_{0})))$.

Note that 
$$
\Cal{V} = \ker(\operatorname{ad}(v_{0}))\bigoplus\ker(\operatorname{ad}(v_{0}))^{\perp}
$$
and 
$$\operatorname{ad}(v_{0}):\, \ker(\operatorname{ad}(v_{0}))^{\perp}\rightarrow 
\operatorname{Im}(\operatorname{ad}(v_{0}))
$$
is an isomorphism.
So, $\operatorname{ad}(v_{0})^{-1}(r_{1}^{-1}(\Cal{Z}_{Q}\cap 
\operatorname{Im}(\operatorname{ad}(v_{0}))))$
is dense in $\ker(ad(v_{0})^{\perp}$.

Denote 
$$A(v_{0}) = 
\operatorname{ad}(v_{0})^{-1}(r_{1}^{-1}(\Cal{Z}_{Q}\cap 
\operatorname{Im}(\operatorname{ad}(v_{0}))).
$$
Fix $u$ such that $r_{1}(v_{0} + [u, v_{0}] + z_{0}) \in \log\Gamma.$
Then, 
$u + A(v_{0}) + \ker(ad(v_{0}))$ is dense in $\Cal{V}$.

For $\bar{u} \in \{u + A(v_{0}) + \ker(ad(v_{0}))\},$ we get,
$$
[\bar{u}, v_{0}] = [u, v_{0}] \mod (r_{1}^{-1}\Cal{Z}_{Q}),
$$
that is, 
$$r_{1}[\bar{u}, v_{0}] = r_{1}[u, v_{0}] \mod (\Cal{Z}_{Q}).$$
Then there is an integer $q$ such that
$$
qr_{1}[\bar{u}, v_{0}] = qr_{1}[u, v_{0}] \mod (\log\Gamma).
$$
Therefore, by \rom{(2.6.6)},
$$
qr_{1}(v_{0} + [\bar{u}, v_{0}] + z_{0}) \in \log\Gamma.
$$
Noting that $n\exp(\Cal{Z}) \subseteq N_0$ whenever $n\in N_0$, it follows from \rom{(2.6.5)} that $N_0$ is dense in $N$.

Now, for $n = \exp( x + w) \in \Cal{N}_{0}, \gamma_{n}(t) 
= n\cdot\exp(v_{0} + z_{0})$ descends to a closed geodesic in $\Gamma\backslash
\Cal{N}$, say of length $l$. This is the longest closed geodesic in the free
homotopy class $[\exp(lv_{0} + l[x, v_{0}] + lz_{0})]$. Thus $F$ sends 
$\gamma_{n}(t)$ to a closed geodesic in $\Phi(\Gamma)\backslash N$ in the
free homotopy class $[\exp(lv_{0} + l[x + \xi, v_{0}] + lz_{0})]$ which
must be the longest one, too.

According to Lemma 1.9, this geodesic must have the form: 
$n_{*}\exp(t(v_{0} + z_{0}))$ for some $n_* \in N$, that is, 
$\tilde{F}(n, v_{0} + z_{0}) = ( n_{*}, v_{0} + z_{0} ).$
Therefore, we have proved that 
$$ 
H(n, v_0 + z_0) = z_0 \quad\text{and}\quad I(n, v_{0} + z_{0}) = v_0,\, \text{for}\, n \in \Cal{N}_{0}.
$$
Since $\Cal{N}_{0} $ is dense in $N$, the above identities hold for all 
$ n \in N.$

Let $ \epsilon \rightarrow 0$ in \rom{(2.6.3)}. Recalling that  
$\|z_{\epsilon, 0}\| < \epsilon$ and $k_{\epsilon} > 0 $ is an integer, we have
$$
\align
\|z_{0}\| &\leq \frac{\|z_{\epsilon, 0}\|}{rk_{\epsilon}} \leq
\frac{z_{\epsilon, 0}}{r} \rightarrow 0. \\
v_{0} &= \frac{rk_{\epsilon}v}{\|rk_{\epsilon}v + z_{\epsilon, 0}\|}
       = \frac{v}{\|v + \frac{z_{\epsilon, 0}}{rk_{\epsilon}}\|} 
\rightarrow v.
\endalign
$$
So, 
$$
I(n, v) = \lim_{\epsilon\rightarrow 0} I(n, v_{0}) = v \,\,
\text{and} \,\, H(n, v) = 0
$$
for all $v\in \Cal{V}_{\Gamma}.$ As noted above, the proposition follows.
\enddemo

\smallskip

\proclaim{2.8 Proposition} For 
$n \in N, v \in \Cal{V}, z \in \Cal{Z},$
$$
\aligned 
A(G^{t}(n, v)) = A(n, v) \qquad A(G^{t}(n, z)) = A(n, z)\\ 
B(G^{t}(n, v)) = B(n, v) \qquad B(G^{t}(n, z)) = B(n, z)
\endaligned
$$
\endproclaim

\smallskip

\demo{Proof}  Immediate from Propositions 2.5 and 2.6.
\enddemo

\smallskip

\proclaim{2.9  Notation and Remarks}
The derived group $[N, N]$ of $N$ is a simply-connected central subgroup with Lie algebra $[\Cal{N}, \Cal{N}] = \Cal{Z}$, and $\exp: \Cal{Z} \rightarrow 
[N, N]$ is a vector space isomorphism. By the remarks in \rom{1.10}, $\Gamma\cap[N, N]\backslash[N, N]$ is a torus; we denote it by $T$. Note that $T$ acts isometrically on $\Gamma\backslash N$ by left translations.

By \rom{Proposition 2.4}, $B(\gamma n, u) = B(n, u), \, 
I(\gamma n, u) = I(n, u)$ for all $\gamma \in \Gamma$, and $A(\gamma n, u) = A(n, u)$ for $\gamma \in [N, N] \cap \Gamma.$ Thus we can define the averages over $T$:
$$ 
\align
\bar{B}(n, u) &= \int_{T}\,
                      B(x\cdot n,\, u)\,dx,  \\
 \bar{A}(n, u) &= \int_{T}\,
                      A(x\cdot n,\, u)\,dx.  \\
 \bar{I}(n, u) &= \int_{T}\,
                       I(x\cdot n, \, u)\, dx.
\endalign
$$
for $(n, u) \in S(N).$ Here $dx$ is the Haar measure of total volume one.
\endproclaim 

\smallskip

\proclaim{2.10 Proposition} In the notation of \rom{1.1, 2.3 } and \rom{2.9},
$$
\align
&\phi(v) = [ \bar{B}(n, \frac{v}{|v|}), \,\,v], \\
&\bar{B}(n, v) = \bar{B}(e, v)
\endalign
$$
for all $(n, v) \in S(N)$ with  $v \in {\Cal V}$.  
\endproclaim

\smallskip

\demo{Proof}  By Proposition 2.4, 
for  all $n \in N, \, v\in \Cal{V},$ with $\|v\| = 1,$ and $ \gamma \in \Gamma,$
we have
$$
\phi(\log \gamma) 
= A(\gamma n, v)\, -\, A(n, v) + [B(n, v),\, \log \gamma],
$$
and therefore,
$$
\phi(\log \gamma) 
= \bar{A}(\gamma n, v)\, -\, \bar{A}(n, v) + [\bar{B}(n, v),\, \log \gamma].
$$
Letting $\pi_v: \Cal{N} \rightarrow \Cal{V}$ be the orthogonal projection and recalling that $\phi(z) = 0$,  we thus have, for $ \eta \in \pi_v(\log(\Gamma))$,
$$
\phi(\eta) = \bar{A}(\exp(\eta)n, v)\, - \, \bar{A}(n, v)
                  +[\bar{B}(n, v), \, \eta].
$$
Next note that $N/[N, N]$ is a simply-connected abelian Lie group. Letting $\pi: N\rightarrow N/[N, N]$ be the projection, then $\pi(\Gamma)\backslash \pi(N)$ is a torus $\Bar{T}$. The Lie algebra of $N/[N, N]$ is $\Cal{N}\backslash\Cal{Z}$.
Under the identification of $\Cal{N}\backslash\Cal{Z}$ with $\Cal{V}$, the exponential map $\Cal{V} \rightarrow N/[N, N]$ carries $\pi_v(\log\Gamma)$ isometrically to $\pi(\Gamma).$ Thus the torus $\Bar{T}$ is isomorphic to $\pi_v(\log\Gamma)\backslash\Cal{V}.$

Since $\bar{B}(yn, v) = \bar{B}(n, v)$ for $y\in[N, N]$ and 
$\bar{B}(\gamma n, v) = \bar{B}(n, v)$ by Proposition 2.4, the map 
$\bar{B}(\cdot, v)$ may be viewed as a map on the torus $\Bar{T}$. Since 
$G^t(n, v) = (n\exp(tv), v)$, Proposition 2.8 implies that 
$\bar{B}(n\exp(tv), v) = \bar{B}(n, v)$  for all $t\in \Bbb{R}$. In particular,
if $v\in \Cal{V}$ is a unit vector such that the projection of $v$ to the torus $\pi_v(\log\Gamma)\backslash\Cal{V}\cong \Bar{T}$ is a generator of the torus 
- i.e., the projection of $\{tv | t \in \Bbb{R}\}$ is dense in the torus - then 
$$
\bar{B}(n, v) \, = \, \bar{B}(n\exp(tv), v)\, \equiv  \, \bar{B}(e, v).
$$
By Kronecker's theorem, the generators form a dense set in the torus. Hence 
$$
\bar{B}(n, v) \equiv  \bar{B}(e, v)
$$
for all unit vectors $v \in \Cal{V}.$
(Note: we can't expect to have  $\bar{A}(n, v) = \bar{A}(e, v)$ since
$\bar{A}(\gamma n, v) \neq \bar{A}(n, v)$ for general $\gamma \in \Gamma$ .)

Hence, for $\eta \in \pi_{v}(\log \Gamma), \, n \in N, 
\, \text{and}\, v \in \Cal{V} \, \text{with} \, \|v\| = 1,$
$$
\phi(\eta) = \bar{A}(\exp(\eta)n, v)\, - \, \bar{A}(n, v)\,
                  +\, [\bar{B}(e, v), \, \eta].
$$
From this identity, we get that for all $n \in N$,
$$
\bar{A}(\exp(\eta)n, v)\, - \, \bar{A}(n, v)\, = \,
\bar{A}(\exp(\eta), v)\, - \, \bar{A}(e, v).
$$
In particular, for $\xi \in \pi_{v}\log(\Gamma)$,
$$
\bar{A}(\exp(\eta)\exp(\xi), v)\, - \, \bar{A}(\exp(\xi)), v)\, = \,
\bar{A}(\exp(\eta), v)\, - \, \bar{A}(e, v).
$$
Therefore, 
$$
\bar{A}(\exp(\eta)\exp(\xi), v)\, - \, \bar{A}(e, v)
= \, (\bar{A}(\exp(\eta), v)\, - \, \bar{A}(e, v))
\, +\, (\bar{A}(\exp(\xi)), v)\, - \,  \bar{A}(e, v)).
$$
Let 
$$
E(\eta, v)\, =\, (\bar{A}(\exp(\eta), v)\, - \, \bar{A}(e, v)),\,\,
\text{for}\, \eta\, \in \pi_{v}\log\Gamma.
$$
Then, 
$$
E(\eta + \xi, v) = E(\eta, v) + E(\xi, v),\,\, 
\text{for}\, \eta\,\,, \xi \, \in \pi_{v}\log(\Gamma).
$$
Let 
$$
\Cal{V}_{Q} = \operatorname{span}_{Q} \{\pi_{v}\log\Gamma \} \subseteq \Cal{V}
$$ 
and let $v_{1}, \cdots, v_{p} \in \pi_{v}\log \Gamma$ be a basis
of $\Cal{V}_{Q}$.

For $w = \Sigma_{i}\, q_{i} v_{i} = \frac{1}{k}\Sigma_{i}\, k_{i} v_{i}
\in \Cal{V}_{\bold{Q}}$, where $k, k_{i}$ are integers, and for all
$v \in \Cal{V}$ with $\|v\| = 1$,
$$
\phi(w) = \frac{1}{k} \phi(\Sigma_{i}k_{i}v_{i})
= \frac{1}{k}E(\Sigma_{i}k_{i}v_{i}, v) + [\bar{B}(e, v), w].
$$
In particular, for $v = \frac{w}{\|w\|} 
= \frac{\Sigma_{i}k_{i}v_{i}}{\|\Sigma_{i}k_{i}v_{i}\|},$ 
since 
$$\bar{A}(\exp(tv), v) = \bar{A}(e, v),
$$ 
we get, 
$$\align
\phi(w) &= \frac{1}{k}(\bar{A}(\exp(\Sigma_{i}k_{i}v_{i}), 
                 \frac{\Sigma_{i}k_{i}v_{i}}{\|\Sigma_{i}k_{i}v_{i}\|})
        - \bar{A}(e, \frac{\Sigma_{i}k_{i}v_{i}}{\|\Sigma_{i}k_{i}v_{i}\|}))
  + [ \bar{B}(e, \frac{w}{\|w\|}), w] \\
   &= [ \bar{B}(e, \frac{w}{\|w\|}), w].
\endalign
$$
Noting that $\Cal{V}_{Q}$ is dense in $\Cal{V}$ and $\bar{B}$ is continuous,
we have proven that for all $v \in \Cal{V}$ with  $v \neq 0$,
$$
\phi(v) = [\bar{B}(e, \frac{v}{\|v\|}),\, v].
$$
\enddemo

\smallskip

Proposition 2.10 gives us a strengthening of Eberlein's result (Proposition 2.2)
as follows:

\proclaim{2.11 Definition} Let $\phi$ be an almost inner derivation of a \rom{2}-step nilpotent Lie algebra $\Cal{N}$. By the definition, there exists a map $\xi: \Cal{N} \rightarrow \Cal{N}$ such that $\phi(x) = [\xi(x), x]$ for all $x\in\Cal{N}$. Note that $\xi$ is not uniquely defined. For $x\in \Cal{Z}$, $\xi(x)$ is completely arbitrary, thus only the values of $\xi$ on $\Cal{V}$ are of interest. We will say $\phi$ is {\it of continuous type} if $\xi$ can be chosen to be continuous on $\Cal{V}\setminus\{0\}$. We will also say the corresponding almost inner automorphism $\Phi$ of $N$ (with differential $\operatorname{\Id } + \phi$ as in  \rom{1.12 (a)}) is of continuous type. \endproclaim

\smallskip

\proclaim{2.12 Theorem}
Suppose $(\Gamma\backslash N, g)$ and  $(\Gamma^*\backslash N^*, g^*)$ are compact \rom{2}-step nilmanifolds whose geodesic flows are $C^0$-conjugate. Then there exists an almost-inner automorphism $\Phi$ of $N$ of continuous type such that $(\Gamma^*\backslash N^*, g^*)$ is isometric to 
$(\Phi(\Gamma)\backslash N, g)$.
\endproclaim

\smallskip

The theorem follows immediately from Propositions 2.2 and 2.10.

\smallskip

\proclaim{2.13 Corollary} There exist \rom{2}-step nilmanifolds $(\Gamma\backslash N, g)$ satisfying:
\roster
\item"(a)" Any \rom{2}-step nilmanifold whose geodesic flow is $C^0$-conjugate to $(\Gamma\backslash N, g)$ must be isometric to $(\Gamma\backslash N, g)$.

\item"(b)" $(\Gamma\backslash N, g)$ is isospectrally deformable, i.e., there exists a continuous family $M_t$ of \rom{2}-step nilmanifolds with 
$M_0 = (\Gamma\backslash N, g)$ such that $M_t$ is isospectral but not isometric to $M_0$ for all $t$.
\endroster
\endproclaim

\smallskip

\demo{Proof}  We only need to find a nilmanifold $(\Gamma\backslash N, g)$
such that (i) every almost inner derivation of $\Cal{N}$ of continuous type is inner and (ii) $\Cal N$ admits a $\Gamma$-almost inner derivation $\phi$ which is non-inner. By Remark 1.14 and Theorem 2.12, condition (i) implies (a). By 1.12, 1.13 and 1.14, condition (ii) implies (b). We will exhibit an example below.
\enddemo

\smallskip

\proclaim{2.14 Example}
 \rom{(see \cite{GM} and \cite{GW}).}  Let $N$ be the six-dimensional simply-connected 
nilpotent Lie group with Lie algebra 
$${\Cal N} 
= \operatorname{span}\{X_{1}, X_{2}, Y_{1}, Y_{2}, Z_{1}, Z_{2}\}
$$
with 
$$
[X_{1}, Y_{1}] = [X_{2}, Y_{2}] = Z_{1}, \quad [X_{1}, Y_{2}] = Z_{2}
$$
and with all other brackets of basis vectors trivial. 

The almost inner
derivations of $\Cal N$  are the linear maps which send   
$X_{2}$ and $Y_{1}$  to multiples 
of $Z_{1}$, send $X_{1}$  and $Y_{2}$  into $\operatorname{span}\{Z_{1}, Z_{2}\}$ and send     
$Z_{1}$ and $Z_{2}$ to $0$. These form a 6-dimensional 
subspace containing the inner derivations as a four-dimensional subspace.
The derivations $\phi_{1}$,   respectively  $\phi_{2}$, 
which send $X_{1}$,    respectively  $Y_{2}$,  to  $Z_{2}$  
and send all other basis vectors to zero span a two-dimensional family of
non-inner almost inner derivations. No almost inner derivation of    
$\Cal N$ is of continuous type. If $M = (\Gamma\backslash N,\, g)$   is any 
nilmanifold associated with $N$, then there
exists a continuous two-parameter family of manifolds strongly Laplace
isospectral 
to $M$ but not $C^{0}$-geodesically conjugate to $M$. 
\endproclaim

See \cite{GM} for further examples.

\smallskip

\smallskip

\heading 3. Special Classes of Nilmanifolds\endheading

\smallskip

In this section, we focus on the special classes of 2-step nilmanifolds introduced in Definition 1.3.  We shall prove the following:

\proclaim{3.1 Theorem} If the geodesic flow of a nilmanifold 
$(\Gamma^{*}\backslash N^{*}, g^*)$ is $C^{0}$-conjugate to that of a nilmanifold
$(\Gamma\backslash N, g)$ of Heisenberg type, then $ (\Gamma^{*}\backslash N^{*}, g^*)$
is isometric to $(\Gamma\backslash N, g)$.
\endproclaim

\smallskip

\proclaim{3.2 Theorem}
If the geodesic flow of a \rom{2}-step nilmanifold $(\Gamma^{*}\backslash N^{*}, g^*)$ is   
$C^{2}$-conjugate to that of a \rom{2}-step nilmanifold $(\Gamma\backslash N, g)$  
which is  strongly in resonance, then $(\Gamma^{*}\backslash N^{*}, g^*)$ is isometric 
to $(\Gamma\backslash N, g)$.
\endproclaim

\smallskip

\proclaim{3.3 Theorem}
If the geodesic flow of a \rom{2}-step nilmanifold $(\Gamma^{*}\backslash N^{*}, g^*)$ is   
$C^{2}$-conjugate to that of an irrational  \rom{2}-step nilmanifold $(\Gamma\backslash N, g)$, then $(\Gamma^{*}\backslash N^{*}, g^*)$ is isometric to $(\Gamma\backslash N, g)$.
\endproclaim

\smallskip

Note that irrationality is a generic condition for 2-step nilmanifolds, so Theorem 3.3 says that generic 2-step nilmanifolds are geodesically rigid.

We continue to use the conventions of 2.3. In particular, we replace 
$\Gamma^{*}\backslash N^{*}$ by $(\Phi(\Gamma)\backslash N, g)$ where $\Phi$ is a $\Gamma$-almost inner automorphism of $N$. To prove Theorems 3.1 - 3.3, we must in each case show that the associated almost inner derivation $\phi$ is inner. The key to the proofs is the following result, valid for arbitrary 2-step
nilmanifolds.

\proclaim{3.4 Proposition}
Let $(\Gamma\backslash N, g)$ be a \rom{2}-step nilmanifold and let $\phi$ be a $\Gamma$-almost inner derivation of the associated Lie algebra $\Cal{N}$. Thus in the notation of \rom{1.1}, there exists a map $\xi: \Cal{V}\rightarrow\Cal{V}$ such that $\phi(x) = [\xi(x),\, x]$ for all $x\in\Cal{V}$. If $\xi$ can be chosen so that $\xi(e^{J(z)}x) = \xi(x)$ for all
$x\in \Cal{V}, \, z\in\Cal{Z},$ then $\phi$ is an inner derivation.
\endproclaim 

\smallskip

\demo{Proof} Pick $v_{1}\in {\Cal V}$ and denote by $V_{1}$ the subspace of 
${\Cal V}$ spanned by all vectors 
$$
\{e^{J(z_{1})}e^{J(z_{2})}\cdots e^{J(z_{k})}v_{1}|\,\,z_{1}, z_{2}, \cdots, 
z_{k} \in {\Cal Z}, k = 0, 1, 2,\cdots \}.
$$ 
Then, pick $v_{2}\in {\Cal V}$ such that $v_{2}\perp V_{1}$. Denote by $V_{2}$
the subspace of ${\Cal V}$ spanned by
$$
\{e^{J(z_{1})}e^{J(z_{2})}\cdots e^{J(z_{k})}v_{2}|\,\,z_{1}, z_{2}, \cdots, 
z_{k} \in {\Cal Z}, k = 0, 1, 2,\cdots \}.
$$ 
We claim that
$$
V_{1} \perp V_{2} \quad \text{and} \quad [V_{1}, V_{2}] = 0. \tag 3.4.1
$$
In fact, let $x = e^{J(z_{1})}\cdots e^{J(z_{k})}v_{1},\,
y = e^{J(\tilde{z}_{1})}\cdots e^{J(\tilde{z}_{l})}v_{2}$. Then

$$
\aligned
<x, \,y> &= <e^{J(z_{1})}\cdots e^{J(z_{k})}v_{1}, \,
e^{J(\tilde{z}_{1})}\cdots e^{J(\tilde{z}_{l})}v_{2}>  \\
        &=< e^{J(-\tilde{z}_{1})}\cdots e^{J(-\tilde{z}_{l})}
        e^{J(z_{1})}\cdots e^{J(z_{k})}v_{1}, \, v_{2}>  \\
        &= 0,
\endaligned
$$
since $v_2 \perp V_1.$

Similarly, we can prove that for any $z\in {\Cal Z}$, 
$<e^{tJ(z)}x,\, y> = 0$. Taking the derivative with respect to $t$ at $t=0$ in this 
identity, we get $<J(z)x,\, y>= 0$, that is $<z,\, [x, y]> = 0.$
Since $z\in {\Cal Z}$ is arbitrary, we conclude that $[x,\, y] = 0$.
Therefore, we have confirmed our claim.

Repeating this process, we can find a finite sequence of vectors 
$v_{1}, v_{2}, \cdots, v_{p}
\in {\Cal V}$ such that
$$
{\Cal V} = V_{1}\bigoplus V_{2}\bigoplus \cdots \bigoplus V_{p}, 
$$
where $V_{i}$ is the subspace spanned by 
$\{e^{J(z_{1})}e^{J(z_{2})}\cdots e^{J(z_{k})}v_{i}|\,\,z_{1}, z_{2}, \cdots, 
z_{k} \in {\Cal Z}, k = 0, 1, 2,\cdots \}$. 
Moreover, $V_{i}\perp V_{j}$ and $[V_{i},\, V_{j}] = 0$ for $i\neq j$.

Since $\phi$ is linear and $\xi(e^{J(z)}v) = \xi(v)$, it is easy to see
that
$$
\phi(x_{i}) = [\xi(v_{i}), x_{i}], \quad \text{for}\,\, x_{i}\in V_{i}.
$$
Let $\xi_{i}$ be the projection of $\xi(v_{i})$ to $V_{i}$. We get,
$$
\phi(x_{i}) = [\xi_{i}, x_{i}], \quad \text{for}\,\, x_{i}\in V_{i}. \tag 3.4.2
$$
Now, for $x\in {\Cal V}$, we write $x = x_{1} + x_{2} + \cdots + x_{p}$, 
where $x_{i} \in V_{i}$. By \rom{3.4.1} and \rom{3.4.2}, 
$$
\aligned
\phi(x) &= \Sigma_{i}\phi(x_{i}) \\
               &= \Sigma_{i}[\xi_{i},\, x_{i}] \\
               &= \Sigma_{i}[\xi, \, x_{i}] \\
               &= [\xi,\, x]     \\
\endaligned
$$
This proves that $\phi$ is  an inner derivation of ${\Cal N}$.
\enddemo

\smallskip

We next prove a technical result valid for arbitrary non-singular 2-step nilmanifolds. This Lemma will be used in the proofs of Theorems 3.2 and 3.3.

\proclaim{3.5 Lemma} We use the notation of \rom{2.1} and \rom{2.3} and assume that $(\Gamma\backslash N, g)$ is a non-singular \rom{2}-step nilmanifold and
$F: S(\Gamma\backslash N) \rightarrow S(\Phi(\Gamma)\backslash N)$ is a $C^2$-geodesic conjugacy. Then for $(n, v+ z)\in S(N)$, with $v\in \Cal{V}, z\in\Cal{Z}$, we have:
\roster
 \item"(a)" $\lim_{s\rightarrow 0}\frac{H(n, \cos(s)v +\sin(s)z)}{s} $

$=\frac{d}{ds}|_{s=0}
 H(n, \cos(s)v +\sin(s)z) = z.$  
\item" "
 \item"(b)" $\lim_{s\rightarrow 0}\frac{d}{ds}
 (\frac{H(n, \cos(s)v +\sin(s)z)}{s}) = 0.$ 
\item" "
 \item"(c)" $ e^{tJ(z)}\frac{d}{ds}|_{s=0}
        I(n, \cos(s)v +\sin(s)z) $
  
$ =\lim_{s\rightarrow 0}\frac{d}{\operatorname{d\lambda}} 
 |_{\lambda=0}I(n\sigma(\frac{t}{s},\cos(s)v +\sin(s)z),\, 
 e^{tJ(z)}\cos(\lambda)v +\sin(\lambda)z).$

\noindent where $\sigma(t, v+z)$ is the geodesic in $N$ with $\sigma(0) = e, \, 
\sigma'(0) = v + z.$
\endroster
\endproclaim

\demo{Proof} (a) The first equality in \rom{(i)} follows from \rom{Proposition 2.6}.
In particular, the limit exists. Write $v(s) = \cos(s) v$ and $z(s) = \sin(s)z$.
 According to Proposition 2.5,  
 $$
 \align
 I&(n\sigma(t, v(s)+z(s)), e^{tJ(z(s))}v(s) + z(s)) \\
 =& 
 e^{tJ(H(n, v(s)+z(s)))}I(n, v(s)+z(s)).  
\endalign
$$

Replacing $t$ by $\frac{t}{s}$ in the above identity, we get
$$
I(n\sigma(\frac{t}{s}, v(s)+z(s)), e^{tJ(\frac{z(s)}{s})}v(s) + z(s)) = 
e^{tJ(\frac{H(n, v(s)+z(s))}{s})}I(n, v(s)+z(s)). \tag 3.5.1
$$
Now, for an arbitrary sequence $s_{n} \rightarrow 0$, there exists a 
subsequence $s_{n_{i}}$ such that $n\sigma(\frac{t}{s_{n_{i}}}, 
v(s_{n_{i}})+z(s_{n_{i}})) \rightarrow g \,\,\operatorname{mod}(\Gamma)$
for some element $g$ in $N$( $g$ may depend on the subsequence).
Therefore, from $I(\gamma n, v+z) = I(n, v+z)$ for $\gamma \in \Gamma$,
we get,
$$
\align
\lim_{i\rightarrow \infty}&
 I(n\sigma(\frac{t}{s_{n_{i}}}, v(s_{n_{i}})+z(s_{n_{i}})), 
 e^{tJ(\frac{z(s_{n_{i}})}{s_{n_{i}}})}v(s_{n_{i}}) + z(s_{n_{i}}))\\
 =&I(g, e^{tJ(z)}v) = e^{tJ(z)}v.
 \endalign
$$
Therefore, 
$$
\lim_{s\rightarrow 0}
 I(n\sigma(\frac{t}{s}, v(s)+z(s)), 
 e^{tJ(\frac{z(s)}{s})}v(s) + z(s))=
 e^{tJ(z)}v.
$$
Hence by (3.5.1) and Proposition 2.6, 
$$ 
\lim_{s\rightarrow 0} e^{tJ(\frac{H(n, v(s)+z(s))}{s})}v=
 \lim_{s\rightarrow 0} e^{tJ(\frac{H(n, v(s)+z(s))}{s})}I(n, v(s)+z(s))
 =   e^{tJ(z)}v,
$$
so, 
$$
 e^{tJ(\lim_{s\rightarrow 0}\frac{H(n, v(s)+z(s))}{s})}v   =   e^{tJ(z)}v.
 $$
Therefore, 
$$
J(\lim_{s\rightarrow 0}\frac{H(n, v(s)+z(s))}{s})v = J(z)v,
$$ 
that is, 
$$
J(\lim_{s\rightarrow 0}\frac{H(n, v(s)+z(s))}{s} - z)v = 0.
$$
Since $N$ is nonsingular,
$ \lim_{s\rightarrow 0}\frac{H(n, v(s)+z(s))}{s} = z.$ 

\smallskip

\noindent (b) First, we show that
$ \lim_{s\rightarrow 0}\frac{d}{ds}
 (\frac{H(n, \cos(s)v +\sin(s)z)}{s}) $  exists.

In fact, by (a), 
$\lim_{s\rightarrow 0} \frac{H(n, \cos(s)v +\sin(s)z)}{s} = z$, so
$$
\align 
\lim_{s\rightarrow 0}&\frac{d}{ds}
 (\frac{H(n, \cos(s)v +\sin(s)z)}{s})\\ 
 &=\lim_{s\rightarrow 0}\frac{1}{s}(\frac{H(n, \cos(s)v +\sin(s)z)}{s}- z) \\
    &=\frac{1}{2}\frac{d^{2}}{ds^{2}}
   |_{s=0}H(n, \cos(s)v+\sin(s)z).  
\endalign
$$
The last equality holds since by Taylor's formula, we have 
$$
\align
&\frac{1}{s}(\frac{H(n, \cos(s)v +\sin(s)z)}{s}- z) \\
&=\frac{1}{2}
   \frac{d^{2}}{ds^{2}}
   H(n, \cos(\bar{s})v+\sin(\bar{s})z) \qquad\text{for some}\,\, 0\leq \bar{s}\leq s.
\endalign
$$

Now, taking the derivative with respect to  $s$ in the identity (3.5.1) and then letting
$s\rightarrow 0$,
we get,
$$
\align
\frac{d}{ds}|_{s=0}& 
I(n\sigma(\frac{t}{s}, v(s)+z(s)), e^{tJ(\frac{z(s)}{s})}v(s) + z(s)) \\
  =&     
 te^{tJ(z)}J(\frac{d}{ds}|_{s=0}
 \frac{H(n, v(s)+z(s))}{s})v \\
 +& 
 e^{tJ(z)}\frac{d}{ds}|_{s=0}I(n, v(s)+z(s)).
\endalign
$$
Equivalently,
$$
\align
J(\frac{d}{ds}|_{s=0}&
 \frac{H(n, v(s)+z(s))}{s})v  \\
 =& - \frac{1}{t}
\frac{d}{ds}|_{s=0}I(n, v(s)+z(s)) \\
+& \frac{1}{t} e^{-tJ(z)}\frac{d}{ds}|_{s=0} 
I(n\sigma(\frac{t}{s}, v(s)+z(s)), e^{tJ(\frac{z(s)}{s})}v(s) + z(s)).
\tag 3.5.2
\endalign
$$   
Now, by Proposition 2.6, 
$$
\align
\frac{d}{ds}|_{s=0}& 
I(n\sigma(\frac{t}{s}, v(s)+z(s)), e^{tJ(\frac{z(s)}{s})}v(s) + z(s)) \\
=&\lim_{s\rightarrow 0}\frac
{I(n\sigma(\frac{t}{s}, v(s)+z(s)), e^{tJ(\frac{z(s)}{s})}v(s) + z(s))
 -e^{tJ(z)}v}{s} \\
=&\lim_{s\rightarrow 0}\frac
{I(n\sigma(\frac{t}{s}, v(s)+z(s)), e^{tJ(\frac{z(s)}{s})}v(s) + z(s))
 - e^{tJ(\frac{z(s)}{s})}v(s)}{s} \\
=&\lim_{s\rightarrow 0}\frac
{I(n\sigma(\frac{t}{s}, v(s)+z(s)), e^{tJ(\frac{z(s)}{s})}v(s) + z(s)) 
 -I(n\sigma(\frac{t}{s}, v(s)+z(s)), e^{tJ(\frac{z(s)}{s})}v(s))}{s} \\
=&\lim_{s\rightarrow 0}
\frac{d}{\operatorname{d\lambda}}|_{\lambda=\bar{s}}
I(n\sigma(\frac{t}{s}, v(s)+z(s)), e^{tJ(\frac{z(s)}{s})}v(s) + z(\lambda)) 
\quad \text{where}\,\, 0\leq\bar{s}\leq s    
\endalign
$$
Using the above equality and the facts that $I(\gamma n, v+z) = I(n, v+z)$ and $F$ is $C^{2}$, we can prove that 
$$
\frac{d}{ds}|_{s=0} 
I(n\sigma(\frac{t}{s}, v(s)+z(s)), e^{tJ(\frac{z(s)}{s})}v(s) + z(s))
$$ 
is bounded for all $t\in \Bbb{R}$. Furthermore, since the $e^{tJ(z)}$ are unitary 
operators, the right hand side of (3.5.2) goes to $0$ when 
$t\rightarrow \infty$. Therefore, the identity (3.5.2) becomes
$$
J(\frac{d}{ds}|_{s=0}
 \frac{H(n, v(s)+z(s))}{s})v = 0
$$
Since $N$ is non-singular, we get 
$\frac{d}{ds}|_{s=0}
 \frac{H(n, v(s)+z(s))}{s} = 0 $.

\smallskip

\noindent (c) Putting 
$
J(\frac{d}{ds}|_{s=0}
 \frac{H(n, v(s)+z(s))}{s})v = 0
$
back into the identity (3.5.2), we get (c).
\enddemo

\smallskip
 
We now turn to Theorem 3.1.

\proclaim{3.6 Lemma}
In the notation of \rom{2.3}, if $N$ is of Heisenberg type, then $H(n, v+z) \equiv z$ for all $(n, v+z)\in SN,$ where $v\in \Cal{V}, z\in \Cal{Z}$.
\endproclaim

\demo{Proof}  First, we take $(n, v+z) \in S(N)$ such that 
$$
\frac{2\pi k}{|z|}(1 + \frac{|v|^{2}}{2|z|^{2}})z 
\in \log(\Gamma\cap [N, N])
$$
for some integer $k$ (take $k$ to be minimum). The set of  such points $(n, v+z)$
are dense in $S(N)$, so, we only need to prove that $H(n, v+z) = z$ for these 
points.

\noindent Let 
$$
\gamma = \exp( \frac{2\pi k}{|z|}(1 + \frac{|v|^{2}}{2|z|^{2}})z)\in 
\Gamma\cap [N, N].
$$ 
Since, by 1.8,  
$$
\align
G^{t}&(n, v+z)   \\
= &(n\exp[(e^{tJ(z)} - \Id )J^{-1}(z)v + 
t(1 + \frac{|v|^{2}}{2|z|^{2}})z - \frac{\sin(t|z|)|v|^{2}}{2|z|^{3}}z],
e^{tJ(z)}v + z),
\endalign
$$
and since
$$
e^{2\pi k J(\frac{z}{|z|})}
 = \operatorname{\Id },
$$
we see that 
$$
G^{\frac{2\pi k}{|z|}}(n, v+z) = \operatorname{dL}_{\gamma}(n, v+z).
$$
So, 
$$
\Tilde{F}G^{\frac{2\pi k}{|z|}} = 
\Tilde{F}\operatorname{dL}_{\gamma}(n, v+z),
$$
and since $\Phi(\gamma) = \gamma$, as $\gamma$ is central, 2.3 implies that 
$$
G^{\frac{2\pi k}{|z|}}\Tilde{F}(n, v+z) = 
\operatorname{dL}_{\gamma}\Tilde{F}(n, v+z).
$$
Now, denote 
$$
\Tilde{F}(n, v+z) = (n_{1}, v_{1} +z_{1})
$$ 
and 
$$
t_{0} = \frac{2\pi k}{|z|}.
$$ 
We have,
$$
\align
G^{t_{0}}\Tilde{F}(n, v+z) &=    \\
(n_{1}\exp[(e^{t_{0}J(z_{1})} - \Id )J&^{-1}(z_{1})v + 
t_{0}(1 + \frac{|v_{1}|^{2}}{2|z_{1}|^{2}})z_{1} - 
\frac{\sin(t_{0}|z_{1}|)|v_{1}|^{2}}{2|z_{1}|^{3}}z_{1}],\,\,
e^{t_{0}J(z_{1})}v_{1} + z_{1}), \\
\text{and} \qquad\qquad\qquad\qquad & \\
\operatorname{dL}_{\gamma}\Tilde{F}(n, v+z) &= (\gamma n_{1}, v_{1} + z_{1}).
\endalign
$$
Comparing   these  identities, we see that 
$$
e^{t_{0}J(z_{1})} = \Id  
$$
and
$$
t_{0}(1 + \frac{|v_{1}|^{2}}{2|z_{1}|^{2}})z_{1} =
t_{0}(1 + \frac{|v|^{2}}{2|z|^{2}})z.
$$
That is, 
$$
e^{t_{0}J(z)} = e^{t_{0}J(z_{1})} = \Id ,
$$
and
$$
z_{1} = \lambda z 
$$
for some positive real number $\lambda$.
Note that $\lambda$ must be  a positive integer since 
$e^{t_{0}J(z)} = e^{\lambda t_{0}J(z)}$.

Now, we have proved that for each $(n, v+z) \in S(N)$, there exists a positive
integer $\lambda$ such that $H(n, v+z) = \lambda z$. Since $S(N)$ is 
connected, this $\lambda$ is independent of $(n, v+z)$. Furthermore,  
$\lambda = 1$ since $H(n, z) = z$ by Proposition 2.5. 
\enddemo

\smallskip

\demo{Proof \rom{of Theorem 3.1}} As in \rom{2.3}, we may assume that $(\Gamma^*\backslash N^*, g^*) = (\Phi(\Gamma)\backslash N, g)$ for some $\Gamma$-almost inner automorphism $\Phi$ of $N$. We only need to show that the corresponding $\Gamma$-almost inner derivation $\phi$ is inner. 

According to  Proposition 2.10,
$$
\phi(v) = [\bar{B}(e, \frac{v}{|v|}), \, v],
\quad \text{for }\,\, 0 \neq v\in {\Cal V}.
$$
On the other hand, since $\phi$ is linear, we get
$$
\aligned
\phi(v) &= -\phi(-v) \\
               &= -[\bar{B}(e, -\frac{v}{|v|}),\, -v] \\
               &= [\bar{B}(e, -\frac{v}{|v|}),\, v]
\endaligned
$$
Therefore, Writing 
$$\Tilde{B}(v)= 
\frac{1}{2}(\bar{B}(e, \frac{v}{|v|})+\bar{B}(e, -\frac{v}{|v|}), 
$$ 
we get
$$
\phi(v) = 
[\Tilde{B}(v), \, v],
\quad \text{for }\,\, 0 \neq v\in {\Cal V}.
$$
According to Proposition 3.4, we only need to show that
$$
\Tilde{B}(e^{J(z)}v) = \Tilde{B}(v)
$$ 
for all $z\in {\Cal Z},\, v\neq 0
\in{\Cal V}.$

By Proposition 2.5, Lemma 1.7 and Lemma 3.6,  we get, for $(n, v+z)\in S(N)$,
$$
B(G^{t}(n, v+z)) = B(n, v+z) + (e^{tJ(z)} - \Id )J^{-1}(z)(I(n, v+z)- v), 
\tag 3.1.1
$$
and
$$
B(G^{t}(n, -v+z)) = B(n, -v+z) + (e^{tJ(z)} - \Id )J^{-1}(z)(I(n, -v+z)+ v). 
\tag 3.1.2
$$
By Proposition 2.5 and Lemma 3.6, 
$$
I(G^{t}(n, v+z)) = e^{tJ(z)}I(n, v+z).
$$  
Taking $t = \frac{\pi}{|z|},$
we get 
$$
I(n\sigma(\frac{\pi}{|z|}, v+z), -v +z) = - I(n, v+z), \tag 3.1.3
$$
where $\sigma(t, v+ z)$ is the geodesic in $N$ with $\sigma(0) = e,\, \sigma'(0) = v + z.$

By (3.1.1) - (3.1.3), we have 
$$ 
\align
&B(G^{t}(n, v+z)) +  B(G^{t}(n\sigma(\frac{\pi}{|z|}, v+z), -v+z)) \\
&= B(n, v+z) + B(n\sigma(\frac{\pi}{|z|}, v+z), -v+z).
\endalign
$$
Hence,
$$ 
\align
\bar{B}(&G^{t}(n, v+z)) +\bar{B}(G^{t}(n\sigma(\frac{\pi}{|z|}, v+z), -v+z))\\
   = &\bar{B}(n, v+z) + \bar{B}(n\sigma(\frac{\pi}{|z|}, v+z),-v+z).
\endalign
$$
Replacing $v$ by $\cos(s)v$, $ z$ by $\sin(s)z$, $t$ by $\frac{t}{\sin(s)}$
in the above identity and recalling that $\bar{B}(n, v)$ is independent of $n$
(Proposition 2.10), 
we get, by letting $s\rightarrow 0$, 
$$ \bar{B}(e, e^{tJ(z)}v) + \bar{B}(e, -e^{tJ(z)}v)
   = \bar{B}(e, v) + \bar{B}(e, -v).
$$
That is $\Tilde{B}(e^{J(z)}v) =\Tilde{B}(v).$ 
\enddemo 

\smallskip

\demo{Proof \rom{of Theorems 3.2 and 3.3}}

As in the proof of Theorem 3.1, we set
$$
\Tilde{B}(v)= 
\frac{1}{2}(\bar{B}(e, \frac{v}{|v|})+\bar{B}(e, -\frac{v}{|v|}).
$$
We need only show that
$\Tilde{B}(e^{J(z)}v) = \Tilde{B}(v)$ for all $z\in {\Cal Z},\, v\neq 0
\in{\Cal V}$.

Note that under the hypotheses of either theorem, $N$ is non-singular. According to Proposition 2.5 and Lemma 1.7, for $(n, v + z)\in S(N),$  
$$
\align
\bar{B}&(G^{\frac{t}{s}}(n, \cos(s)v +\sin(s)z) -
 \bar{B}(n, \cos(s)v +\sin(s)z) \\
 =&
 (e^{tJ(\frac{H(n, \cos(s)v +\sin(s)z)}{s})}-\Id )
 J^{-1}(H(n, \cos(s)v +\sin(s)z))
 \bar{I}(n, \cos(s)v +\sin(s)z) \\
 -&(e^{tJ(\frac{\sin(s)}{s}z)} - \Id )J^{-1}(\sin(s)z)\cos(s)v \\
 =& ( e^{tJ(\frac{H(n, \cos(s)v + \sin(s)z)}{s})} - \operatorname{\Id })
  J^{-1}(\frac{H(n, \cos(s)v + \sin(s) z}{s}) \frac{\Bar{I}(n, \cos(s)v + \sin(s)z)}{s} \\
-&(e^{tJ(\frac{\sin(s)}{s}z)} - \operatorname{\Id }) J^{-1}(\frac{\sin(s)}{s} z) \frac{\cos(s)}{s} v
\endalign
$$
where $\bar{I}(n, v+z)$ is defined in 2.9. The last equality follows from the fact that $J^{-1}(\frac{z}{s}) = J(\frac{z}{s})^{-1}$.

Letting $s\rightarrow 0$ in the above identity and using Lemma 3.5 (a) and (b)
and Lemma 1.7, 
we get
$$
 \bar{B}(n, e^{tJ(z)}v) - \bar{B}(n, v) =
 (e^{tJ(z)}-\Id )J^{-1}(z)\frac{d}{ds}|_{s=0}
 \bar{I}(n,\cos(s)v +\sin(s)z). \tag 3.2.1
 $$
 By Proposition 2.10,
 $\bar{B}(n, v)$ is independent of $n$, thus so is 
$$
 (e^{tJ(z)}-\Id )J^{-1}(z)\frac{d}{ds}|_{s=0}
 \bar{I}(n,\cos(s)v +\sin(s)z).
$$ 
Therefore, by Lemma 3.5 (c), 
 $$
\align
 e^{tJ(z)}\frac{d}{ds}&|_{s=0}
 \bar{I}(n, \cos(s)v +\sin(s)z)   \\
 = 
\frac{d}{ds}& 
|_{s=0}\bar{I}(n, e^{tJ(z)}\cos(s)v +\sin(s)z).
\endalign  
$$
If $N$ is  strongly in resonance, we can pick a $t$ such that 
$ e^{tJ(z)}=-\Id .$ If $N$ is irrational, there is a sequence $\{t_{i}\}$
such that $\lim_{i\rightarrow \infty}e^{t_{i}J(z)} = -\Id $.
In both cases, we have,
 $$
\align
\frac{d}{ds}& 
|_{s=0}\bar{I}(n, -\cos(s)v +\sin(s)z) \\
= &-\frac{d}{ds}|_{s=0}
 \bar{I}(n, \cos(s)v +\sin(s)z). 
\endalign 
$$
Substituting this  identity into (3.2.1), we get  
$\Tilde{B}(e^{J(z)}v) = \Tilde{B}(v)$. 
This 
completes the proofs of Theorems 3.2 and 3.3.
\enddemo

\enddocument